\documentclass{amsart}

\usepackage[hyperref]{hyperref}

\usepackage{rotating}
\usepackage[latin1]{inputenc}
\usepackage{graphicx}
\usepackage{amsfonts}
\usepackage{amssymb}
\usepackage{amsmath}
\usepackage{bbm}

\numberwithin{equation}{section}

\renewcommand{\phi}{\varphi}
\newcommand{\C}{{\mathbb{C}}}
\newcommand{\R}{{\mathbb{R}}}
\newcommand{\Q}{{\mathbb{Q}}}
\newcommand{\Z}{{\mathbb{Z}}}

\renewcommand{\epsilon}{\varepsilon}
\renewcommand{\theta}{\vartheta}

\DeclareMathOperator{\im}{im}

\newtheorem{theorem}{Theorem}[section]
\newtheorem{lemma}[theorem]{Lemma}
\theoremstyle{definition}

\newtheorem{remark}[theorem]{Remark}

\title[Open books]{Open books on contact five-manifolds}

\author{Otto van Koert}
\address{Université Libre de Bruxelles\\
Département de Mathématiques - CP 218\\
Boulevard du Triomphe\\
B-1050 Bruxelles - Belgique}
\email{ovkoert@ulb.ac.be}

\thanks{This paper is a part of my thesis, which I wrote under supervision of H.~Geiges at the university of Cologne. I am indebted to him for his support and patience. I would also like to thank the anonymous referee for pointing out a nicer argument for a part of the construction. Currently I am supported by the F.N.R.S., Belgium.}



\begin{document}

\begin{abstract}
The aim of this paper is to give an alternative proof of a theorem about the existence of contact structures on five-manifolds due to Geiges. This theorem asserts that simply-connected five-mani-folds admit a contact structure in every homotopy class of almost contact structures. Our proof uses the open book construction of Giroux.
\end{abstract}

\maketitle

\section{Introduction}
\label{introduction}
At the ICM of 2002 Giroux announced his results on the relation between contact manifolds and open book decompositions. The easy part of his results (and the part that we shall use) is a generalization of a construction due to Thurston and Winkelnkemper \cite{Thurston_Winkelnkemper}; one can adapt certain open book decompositions to contact structures, thus giving a procedure to construct contact structures using open books. Roughly speaking Giroux's construction goes as follows. Take a compact Stein manifold $P$ or more generally an exact symplectic manifold with boundary and a symplectomorphism $\psi$ of $P$ that is the identity near the boundary of $P$. The mapping torus of $(P,\psi)$ can be shown to admit a natural contact structure. On the other hand, a neighborhood of the binding $\partial P\times D^2$ has a natural contact structure that can be glued to the contact structure on the mapping torus, therefore giving rise to a closed contact manifold with an adapted contact structure.

In this paper, we will use Giroux's construction to reprove a theorem on the existence of contact structures on five-manifolds due to Geiges \cite{Geiges_five}. More precisely, we shall reprove the following theorem.
\begin{theorem}[Geiges] 
\label{thm_main}
Let $M$ be a simply-connected five-manifold. Then
  $M$ admits a contact structure in every homotopy class of almost contact structures.
\end{theorem}
The main idea of our alternative proof is very simple. Using the classification of simply-connected five-manifolds, we can reduce the problem to finding contact structures on certain model manifolds. We do this by explicit construction using Giroux's procedure. Although this is not necessary in the construction of Giroux, we will always take Stein surfaces as pages. Since the classification of simply-connected five-manifolds is determined by the homology groups and the second Stiefel-Whitney class, it suffices to track these topological invariants. 

\section{Preliminaries}
\label{sec_prelims}
We start by recalling Giroux's construction in a bit more detail. Let $P$ be a compact Stein manifold of real dimension $2n$ and take a strictly plurisubharmonic function $f$. The function $f$ defines an exact symplectic form $d\beta=-d(d^c f)=-d(df\circ J)$, where $J$ is the complex structure on $P$.  Let now $\psi:~P\to P$ be a symplectomorphism that is the identity near the boundary of $P$. In general, $\psi$ does not preserve $\beta$, which we would like to have. However, it turns out that the pull-back of $\beta$ under $\psi$ can be assumed to be exact due to the following lemma of Giroux \cite{Giroux_talk}.
\begin{lemma}[Giroux]
The symplectomorphism $\psi$ can be isotoped to a symplectomorphism $\psi'$ that is the identity near the boundary and that satisfies
$$
\psi'^*\beta=\beta-dh.
$$
\end{lemma}
\begin{proof}
Let us denote the one-form $\psi^*\beta-\beta$ by $\mu$. Since $d\beta$ is
non-degenerate, we find a unique solution $Y$ to the equation
$i_Y d\beta=-\mu$. The flow of the vector field $Y$ preserves $d\beta$,
because $\mu$ is closed,
$$
0=-d\mu=d i_Y d\beta =\mathcal L_Y \beta.
$$
Since $\psi$ is the identity near the boundary, $\mu$ and hence $Y$ vanish
near the boundary. If we denote the time $t$ flow of $Y$ by $\phi_t$, then we
see that $\psi'=\psi\circ\phi_1$ is a symplectomorphism that is the
identity near the boundary. Note that $\mathcal L_Y \mu =0$, so
$\phi_t^*\mu=\mu$ for all $t$. We check that the difference of the pullback
of $\beta$ and $\beta$ is indeed exact. We have
$$
(\psi \circ \phi_1)^* \beta -\beta=\phi_1^*(\mu+\beta)-\beta=\mu+\phi_1^*\beta -\beta.
$$
On the other hand, we can express the difference $\phi_1^* \beta-\beta$ as
\begin{align*}
\phi_1^*\beta-\beta&=\int_0^1\frac{d}{dt} \phi_t^*\beta=\int_0^1 \phi_t^*
\mathcal L_Y \beta=\int_0^1 \phi_t^* \left( i_Y d\beta +d (i_Y\beta)
\right) \\
&=-\mu+\int_0^1
d\phi_t^*(i_Y\beta).
\end{align*}
Moving $\mu$ to the left-hand-side, we see that $\mu+\phi_1^*
\beta -\beta$ is exact, which shows the claim of the lemma. 
\end{proof}
Using this lemma we can make a mapping torus with a natural contact structure. The form
$$
\alpha=d\phi+\beta
$$
is a contact form on $P\times \R$ that descends to the perturbed mapping torus
$$
A:=P\times \R /(x,\phi)\sim \Psi(x,\phi)=(\psi(x),\phi+h(x)).
$$
We see that $\alpha$ indeed gives $A$ a well defined contact form, because
$$
\Psi^*\alpha=d\phi+d h+\psi^*\beta=d\phi+dh+\beta-dh=\alpha.
$$
The boundary of the page $K=\partial P$ inherits a natural contact form $\gamma=\beta|_{TK}$, since $P$ is a compact Stein manifold. We use this to "complete" $A$ into an open book. Glue $B:=K\times D^2$ along its boundary to $A$. This can be done in a natural way, since $\psi$ was assumed to be the identity near the boundary of $P$. 

\begin{figure}
\begin{center}
\include{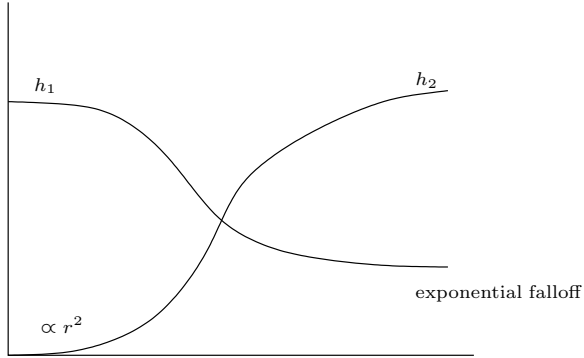}
\end{center}
\caption{The functions $h_1$ and $h_2$}
\label{fig_functions}
\end{figure}

This construction involving a mapping torus is sometimes called an {\bf abstract open book}. Note that one can put a contact form $\tilde \alpha$ on $B$ that matches the contact form on $A$, thus giving rise to a closed contact manifold $X:=A\cup_\partial B$. This contact form $\tilde \alpha$ has the form
$$
\tilde \alpha=h_1(r) \gamma+h_2(r) d\phi,
$$
where $(r,\phi)$ are polar coordinates on $D^2$ and $h_1$ and $h_2$ are functions that are sketched in Figure~\ref{fig_functions}. For the choice of functions indicated in Figure~\ref{fig_functions} the form $\tilde \alpha$ is in fact a contact form, since the contact condition can be rewritten as
$$
\tilde \alpha \wedge d\tilde \alpha^n= h_1^{n-1} \frac{h_1 h_2' -h_2h_1'}{r}\gamma\wedge d\gamma^{n-1}\wedge dr \wedge r d\theta.
$$
This is a non-vanishing form, since $\frac{h_1 h_2' -h_2h_1'}{r}\neq 0$ by our choice of functions $h_1$ and $h_2$. Also note that by choosing these functions suitably, we can ensure that the contact form $\tilde \alpha$ matches the contact form $\alpha$ near the boundary of $A$. Hence we get a well defined contact form on the entire abstract open book. We will call the abstract open book together with the contact form given by the above construction an {\bf abstract contact open book}. In this procedure the contact structure is determined by the page $P$ and the monodromy $\psi$ up to isotopy.

\begin{remark}
\label{rmk_booksum}
If two manifolds, say $M$ and $N$, are constructed via this procedure, then their connected sum $M\# N$ can also be constructed this way. Indeed, if $M$ is an abstract contact open book coming from the pair $(P_1,\psi_1)$ and $N$ is constructed from $(P_2,\psi_2)$, then we may consider the boundary connected sum $P_1 \natural P_2$, which is again a Stein manifold. Note that the symplectomorphisms $\psi_1$ and $\psi_2$ can be glued to a symplectomorphism $\psi_1\natural \psi_2$ of $P_1 \natural P_2$, since both symplectomorphisms are the identity near the boundary. Then the abstract contact open book constructed from $(P_1 \natural P_2,\psi_1\natural \psi_2)$ provides an open book decomposition for $M\# N$. This procedure is called a {\bf book-connected sum}.
\end{remark}

\subsection{Classification of simply-connected five-manifolds}
\label{sec_classification}
We now recall Barden's classification of simply-connected five-manifolds \cite{Barden}. For a simply-connected manifold $M$ we can regard the second Stiefel-Whitney class as a map $w_2(M):~H_2(M)\to \Z_2$.
\begin{theorem}[Barden] Two simply-connected five-manifolds $M_1$ and $M_2$ are diffeomorphic if and only if there exists an isomorphism of groups $A:~H_2(M_1)\to H_2(M_2)$ such that
$$
w_2(M_1)=w_2(M_2)\circ A.
$$
\end{theorem}
Before we give a description of the decomposition of a simply-connected five-manifold into prime manifolds, we would like to point out that a necessary condition for the existence of a contact form is the existence of an almost contact structure.  The existence of an almost contact structure is governed by purely topological considerations. For instance, a simply-connected five-manifold $M$ admits an almost contact structure if and only if the third integral Stiefel-Whitney class $W_3(M)=0$, see Lemma 7 from \cite{Geiges_five}. 

A simply-connected five-manifold can be uniquely decomposed into a connected sum of prime manifolds $M_k$ for $1 \leq k \leq \infty$ with possibly one extra summand $X_j$ with $j=-1$ or $1 \leq j \leq \infty$. The second Stiefel-Whitney class of $X_j$, the class $w_2(X_j)$, is always non-trivial.

The manifold $M_k$ has homology group $H_2(M_k)\cong \Z_k \oplus \Z_k$ for $1<k<\infty$. The manifold $M_\infty$ can be identified with $S^2\times S^3$. In the decomposition above we always take $k$ to be a prime power if $k\neq \infty$. The manifold $M_1$ is $S^5$ and is only needed in a decomposition of $M$ if $M\cong S^5$. These manifolds all carry an almost contact structure since $W_3(M_k)=0$. 

The manifold $X_{-1}$ is known as the Wu-manifold and satisfies $H_2(X_{-1})=\Z_2$. It does not carry an almost contact structure since $W_3(X_{-1})\neq 0$. For $1\leq j <\infty$ we have $H_2(X_j)=\Z_{2^j}\oplus \Z_{2^j}$. Again $W_3(X_j)\neq 0$, so we do not need to consider these manifolds because they cannot have a contact structure. Finally the manifold $X_\infty$ can be identified with $S^2\tilde \times S^3$, the non-trivial $S^3-$bundle over $S^2$ and has $H_2(X_\infty)\cong \Z$. Among the "$X$"-manifolds $X_\infty$ is the only one with vanishing $W_3$, so we shall need to consider $S^2\tilde \times S^3$.

Using this decomposition we see that it suffices to compute the second homology group and the second Stiefel-Whitney class in order to determine which contact five-manifold we have.

\subsection{Some general arguments for computing the homology of open books}
\label{sec_general}
In our construction we will always use a simply-connected page. This implies that the abstract open book will also be a simply-connected manifold. Indeed, if we use $P$ to denote the page of the open book and $A$ to denote the mapping torus of $P$, we see that the homotopy exact sequence of a fibration implies that $\pi_1(A)=\Z$. Now consider the completed open book $X$, obtained by gluing $B:=\partial P\times D^2$ to $A$ along a collar neighborhood of its boundary. Note that the generator of the fundamental group of $A$ gets killed in $B$; the curve $\{ point\}\times S^1$ lying in the boundary of $A$ represents the generator. In $B$, this curve bounds the disk $\{ point\} \times D^2$. On the other hand, we can always choose a curve lying in $\partial P\times \{ point \}$ to represent a generator of $\pi_1(B)$. However, such a curve will always be contractible in $A$, since it lies in a page. An application of the Seifert-Van Kampen theorem shows that the open book is simply-connected.

Since the classification of simply-connected five-manifolds is mainly controlled by homology, some general arguments to compute the homology of open books turn out to be useful. First of all, we shall stick to the notation introduced in Section \ref{sec_prelims}, namely we shall denote the mapping torus of a compact Stein manifold $P$ by $A$, the thickened binding by $B$ and the closed manifold by $X:=A\cup_\partial B$. We can, in fact, glue along a collar neighborhood of the boundary. Therefore, we can apply the Mayer-Vietoris sequence straight away to $X$ and its "parts" $A$ and $B$ to compute the homology of $X$.

The homology of the mapping torus $A$, being a fiber bundle over $S^1$, can be determined from the Wang sequence \cite{Wang}, but see also \cite{ACampo}. This  works as follows. Suppose $P$ is a manifold and $\phi$ a diffeomorphism of $P$. If the mapping torus $A$ is defined by
$$
A:=P\times [0,1]/{(x,0)\sim(\phi(x),1)},
$$
then we have the following long exact sequence in homology, called the Wang sequence,
$$
\to H_3(A;\Z) \to H_2(P;\Z) \stackrel{\phi_*-id}{\to} H_2(P;\Z) \stackrel{incl_*}{\to} H_2(A;\Z) \to.
$$
The homology of $B$ is simply the homology of the boundary of a page $K=\partial P$. Finally we have the homotopy equivalence $A\cap B \sim K\times S^1$, so the homology of $A\cap B$ can be determined using the K\"unneth formula for $K\times S^1$.

In order to simplify the sequences, we will use the following simple argument. If $\phi:~G \to G$ is a surjective homomorphism of finitely generated abelian groups, then $\phi$ is an isomorphism. This can be seen as follows. Write $G=\Z^k\oplus T$, where $\Z^k$ is a free abelian group of rank $k$ and $T$ is a torsion group. Write $\phi=(f,g)$, where $f:~G\to \Z^k$ and $g:~G\to T$. Of course, $f$ cannot depend on the torsion part of $G$, so $f$ can be regarded as a surjective homomorphism from $\Z^k$ to $\Z^k$. This means $f$ must be injective, since this would also be true if we extended $f$ to a linear surjection from $\Q^k$ to $\Q^k$. This implies that if we restrict $g$ to $T$, we get a surjective map from $T$ to $T$. Since these are finite sets with an equal number of points, the map $g|_T$ must be injective as well, which in turn implies that $\phi$ is injective.

We will apply this for instance in the following situation. Consider the Mayer-Vietoris sequence of the pair $(A,B)$ in $X$, where $A$ and $B$ are as above. Since we already saw that $X$ is simply-connected, we also have that $H_1(X)=0$, and hence a part of the Mayer-Vietoris sequence looks like
$$
H_1(A\cap B) \stackrel{f}{\to} H_1(A)\oplus H_1(B) \to 0. 
$$
Note that by the K\"unneth formula $H_1(A\cap B)\cong H_1(A) \oplus H_1(B)$. Applying the above argument at this point shows that the map $f$ is an isomorphism. This can also be seen in different ways, for instance using the fundamental groups of the involved spaces.

\section{Contact open books for $S^2\times S^3$ and $S^2\tilde \times S^3$}
\label{sec_openbook_free}
\begin{figure}
\begin{center}
\include{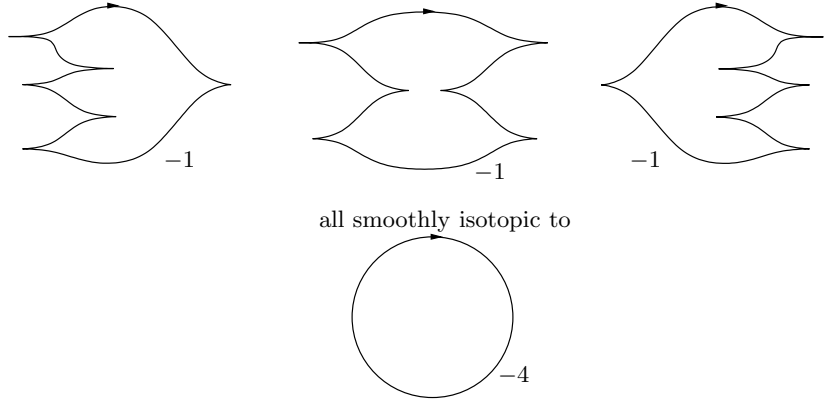}
\end{center}
\caption{Different Stein structures on $\Sigma_4$}
\label{Figure_sigma_k}
\end{figure}

Our construction starts by taking a simple Stein
manifold $P:=\Sigma_k$, the 2-disk-bundle over $S^2$ with Euler number $-k$
with $k\geq 2$. We
remark that these manifolds carry often more than one Stein structure as
can be seen in Figure \ref{Figure_sigma_k}. Here we use the Kirby diagram description of Stein surfaces due to Gompf \cite{Gompf_Stein}; by attaching two-handles in a suitable way to Legendrian knots, one can ensure that the resulting manifold carries a Stein structure, i.e.~we choose the framing of a Legendrian knot $K$ to be equal to the contact framing minus $1$.
First we will show that we get contact open books for $S^2\times S^3$ and
$S^2\tilde \times S^3$, then we will show that the different realizations
from Figure \ref{Figure_sigma_k} can give rise to different contact structures on
$S^2\times S^3$ and $S^2\tilde \times S^3$. 

Let $S_k$ denote the contact boundary of $\Sigma_k$. It is well known that the manifold $S_k$ can be identified with the circle bundle over $S^2$ with Euler
number $-k$. We will use the identity
as monodromy, so the mapping torus of the pair $(\Sigma_k,id)$ is
diffeomorphic to $A:=\Sigma_k\times S^1$. A neighborhood of the binding will be
written as $B:=S_k\times D^2$. By gluing $A$ and $B$ in a collar neighborhood of their boundary we obtain a contact manifold $X:=A\cup_\partial B$.

To see what manifold $X$ is, consider the rank $4$ disk bundle $\Sigma_k\times D^2$ over $S^2$. We can rewrite its boundary as
$$
\partial(\Sigma_k\times D^2)= \Sigma_k \times S^1 \cup_\partial S_k \times D^2=A\cup_\partial B=X.
$$
In other words, the manifold $X$ is a $3$-sphere bundle over $S^2$. To see what sphere bundle it is, we look more closely at the vector bundle associated to the disk bundle $\Sigma_k$, which we shall denote by $\sigma_k$. If we denote the trivial bundle of rank $2$ by $\epsilon^2$, then $\Sigma_k\times D^2$ is the disk bundle associated to $\sigma_k\oplus \epsilon^2$. Recall now that rank $4$ vector bundles over $S^2$ are classified by their second Stiefel-Whitney class. In our case, this class is given by
$$
w_2(\sigma_k\oplus \epsilon^2)=w_2(\sigma_k)=k\mod 2.
$$
So for $k$ even the bundle $\sigma_k\oplus \epsilon^2$ is trivial and for $k$ odd the bundle $\sigma_k\oplus \epsilon^2$ is the unique non-trivial bundle of rank $4$ over $S^2$. As a result, we see that $X$ is diffeomorphic to $S^2\times S^3$ for $k$ even. For $k$ odd, the manifold $X$ is diffeomorphic to $S^2 \tilde \times S^3\cong X_\infty$. 

\subsection{Chern classes of contact structures}
Let us take a look at Figure \ref{Figure_sigma_k}. Legendrian unknots
representing $\Sigma_k$ have rotation numbers going from
$-k+2,-k+4,\cdots,k-2$. Fix a Legendrian unknot representing $\Sigma_k$ and
denote its rotation number by $r$. Now Theorem 11.3.1 from the book of Gompf and Stipsicz \cite{Gompf_Stipsicz} 
tells us how to compute the Chern class.
\begin{theorem}[Gompf]
Suppose $P$ is a Stein surface obtained by two-handle attachment along a Legendrian link $L$. Then $c_1(P)$ is represented by a cocycle whose value on each oriented two-handle $h$ attached along a component $K$ of $L$ is given by $r(K)$.
\end{theorem}
We have just a single Legendrian unknot, so application of this theorem shows that 
$$
c_1(\Sigma_k)=r\in \Z\cong H^2(\Sigma_k).
$$

We now want to establish the relation between the Chern class of the
contact structure corresponding to the open book decomposition we described
and the Chern class of $\Sigma_k$, the page of the open book. We may regard
the pull-back $p_1^*T\Sigma_k$ as a subbundle of $TA$. If we denote the
symplectic form on $\Sigma_k$ by $\omega$, then we may write the contact
form on $A$ as $\alpha=dt+\beta$, where $t$ is the local coordinate on
$S^1=\R/\Z$, and $\beta$ satisfies $d\beta=p_1^*\omega$. We obtain a complex
structure $J$ for $p_1^*T\Sigma_k$ by pulling back the (almost) complex
structure on $\Sigma_k$ that is compatible with $\omega$.

Next, we construct a vector bundle isomorphism from $p_1^*T\Sigma_k$ to the
contact structure $\xi=\ker \alpha$. Define
\begin{eqnarray*}
\phi:~ p_1^*T\Sigma_k & \to & \xi \\
v & \mapsto & v-\beta(v)\frac{\partial}{\partial t}.
\end{eqnarray*}
In the definition of this map, we regard both $p_1^*T\Sigma_k$ and $\xi$ as subbundles of the tangent
bundle. The vector field
$\frac{\partial}{\partial t}$ generates the standard rotation in the
$S^1$-direction. 

The inverse of $\phi$ can be obtained as follows,
$$
\phi^{-1}(v)=H(Tp_1(v)),
$$
where we use $H$ to denote the obvious lift from $T\Sigma_k$ to $TA$. In other
words, the inverse of $\phi$ projects out the $\frac{\partial}{\partial
  t}$-component of an element in $\xi\subset TA$. This map $\phi$ can be
used to give $\xi$ a complex structure. Put $\tilde J:=\phi \circ J \circ \phi^{-1}$.
This makes $\phi$ into complex vector bundle isomorphism from
$(p_1^*T\Sigma_k,J)$ to $(\xi,\tilde J)$, because by
construction $\tilde J \circ \phi=\phi \circ J$. We check now that $\tilde J$ is a complex structure for $\xi$ compatible with
$d\alpha=d\beta$. We set $\tilde v=\phi (v)$ and $\tilde w=\phi (w)$. Then
\begin{align*}
d\beta( \tilde J \tilde v,\tilde J \tilde w)&=d\beta(\phi(J v),\phi
( J w) )=d\beta(Jv,Jw)=d\beta(v,w)\\
&=d\beta( \phi(v), \phi(w) )=d \beta(\tilde v,\tilde w)
\end{align*}
These steps hold true, because $\phi$ adds an $S^1$-component
and $d\beta$ does not contain any $dt$ part, so $d\beta(\phi (\ldots), \phi
(\ldots) )=d\beta (\ldots,\ldots)$. Also, $J$ is a complex structure on
$(p_1^*T\Sigma_k,J)$ compatible with $d\beta$. For the same reasons, the
following holds: 
$$
d\beta(\tilde v,\tilde J \tilde v)= d\beta(\phi(v),\phi(Jv))=d\beta(v,J
v) >0 \text{ if } \tilde v\neq 0.
$$
This proves that $\tilde J$ is a complex structure
compatible with the contact structure $\xi$. Since $(p_1^*T\Sigma_k,J)$ and
$(\xi,\tilde J)$ are isomorphic as complex vector bundles by $\phi$ (which
covers the identity), their Chern classes are the same. We had already
computed the Chern class of $\Sigma_k$, so we have proved that
$c_1(\xi)=r\in \Z\cong H^2(A)$.

We resort to a Mayer-Vietoris argument to complete our computation of
the Chern class of $X$. Consider the Mayer-Vietoris sequence for cohomology
with integer coefficients. The part that is relevant to us looks like
$$
0 \to \underset{\cong \Z}{H^1(A)}\oplus \underset{\cong 0}{H^1(B)} \stackrel{\alpha}{\to}
\underset{\cong \Z}{H^1(A\cap B)} \stackrel{f}{\to} \underset{\cong
  \Z}{H^2(X)} \stackrel{(i^*,j^*)}{\to} \underset{\cong \Z\oplus \Z_k}{H^2(A) \oplus H^2(B)}.
$$
Since the map $\alpha$ is injective, it has to map $1$ to some non-zero
integer, say $m$. If $m$ is not equal to $\pm 1$, then we see that $f(m)=0$,
but $f(1)\neq 0$ by exactness. However $H^2(X)$ has no torsion, so we see
that $m=\pm 1$ and thus the map $\alpha$ is an isomorphism. Again, by
exactness the map $f$ has to be the zero
homomorphism. So we see that the map $(i^*,j^*)$ is injective. We can say a
bit more, namely that $i^*$ is injective. This can be seen
by noting that $H^2(B)$ is torsion. We
show that $i^*$ is an isomorphism by looking at the sequence of the pair
$(X,A)$. The piece of the sequence that interests us, looks like
$$
H^2(X) \stackrel{i^*}{\to} H^2(A) \to H^3(X,A).
$$
By excision, we have $H^3(X,A)\cong H^3(B,\partial B)$. The latter group is
seen to be isomorphic to $H_2(B)= 0$ by Poincar\'e duality. This shows
that $i^*$ is surjective.

The restriction of the first Chern class of the contact structure $\xi_X$ on
$X$ to $A$ is given by $c_1(\xi_X)=r$. Since we just checked $i^*$ to be an
isomorphism, it follows that $c_1(\xi_X)=r\in \Z\cong H^2(X)$. There is an ambiguity in this notation, namely it depends on which
generator of $H^2(X)$ we take.

These ambiguities do not matter for the point we want to make, which is
showing that all possible Chern classes of $\xi_X$ can be realized by our
open books (i.e.~both positive and negative elements in $H^2(X)$). Indeed, the isomorphism $i^*:~H^2(X)\to H^2(A)$ only depends on the topological structure of $\Sigma_k$ and $S_k$,  and not on the Stein structure of $\Sigma_k$. Hence we can change the sign of the first Chern class of $\xi_X$ without affecting the orientation of $X$, for example by replacing the Legendrian knot representing $\Sigma_k$ by its mirror.

Notice that for $(X\cong S^2 \times S^3,\xi_X)$ we can realize all even
Chern classes and for $(X\cong S^2 \tilde \times S^3,\xi_X)$ we can realize
all odd Chern classes. Namely, observe that the rotation number $r$ of the
diagram in Figure \ref{Figure_sigma_k} can attain any even value, provided
that we have chosen $k$ even and large enough for that purpose. The same
argument works for odd rotation numbers.

\section{Open books for prime manifolds}
\label{sec_openbook_torsion}
In this section we will construct open book decompositions of the remaining
prime manifolds, i.e.~those simply-connected five-manifolds with torsion
$H_2$ and trivial Stiefel-Whitney class. In order to cover these remaining cases, we turn our attention to a well studied class of Stein manifolds, namely Brieskorn varieties. Note that for the cases we still need to cover, it is necessary to use a non-trivial monodromy.

\subsection{Brieskorn varieties}
Consider the polynomial 
$$
P_t(z)=\sum_{i=0}^{n}z_i^{a_i}-t
$$
for $z=(z_0,\ldots,z_n)\in \C^{n+1}$ and $t\in \C$. The zero set of this
polynomial is a Stein manifold if $t\neq 0$. If $t=0$, the zero set of
$P_t$ has a singularity at 0 if one of the exponents is larger than 1. We
will denote the zero set of the polynomial $P_t$ by $\Sigma_a$, where $a$
indicates that this set depends on the exponents
$a=(a_0,\ldots,a_n)$. We will call the set $\Sigma_a$ a {\bf Brieskorn variety}. There is a group action of $\Z_{a_i}$ on $\Sigma_a$
obtained by multiplying the $i^\text{th}$ coordinate by $a_i^\text{th}$
roots of unity for each $i=0,\ldots,n$. These Stein manifolds can be made
into compact Stein manifolds by restricting $\Sigma_a$ to a ball $B_R=\{ z\in
\C^{n+1}~|~ |z|\leq R \}$ in $\C^{n+1}$ with sufficiently large radius. By
abuse of notation, we will also denote this set by $\Sigma_a$. The boundary of
this compact Stein manifold is a {\bf Brieskorn manifold} with exponents $a$,
provided that $t$ is small enough. This
property of Brieskorn manifolds can for instance be found in theorem 14.3 of
\cite{Hirzebruch}.

We would like to use Brieskorn varieties as pages with the corresponding
Brieskorn manifolds as binding in open books. In order to produce a non-trivial symplectomorphism, we consider the action of the generator of $\Z_{a_0}$ on
$\Sigma_a$ as monodromy, i.e.~we use the ``rotation'' map
\begin{eqnarray*}
\phi: \Sigma_a  & \to & \Sigma_a \\
(z_0,\ldots,z_n) & \mapsto & (\zeta_{a_0}z_0,z_1,\ldots,z_n),
\end{eqnarray*}
where $\zeta_{a_0}$ is the $a_0^\text{th}$ root of unity $e^{2 \pi i
  /a_0}$. Since this is even a biholomorphism, we
get a symplectomorphism of the page, but we still need to show that we can
isotope this map symplectically to the identity near the boundary of the
  page. We will describe this in the following interlude.

\subsubsection{The rotation maps $\phi$ are symplectically isotopic to the
  identity}
\label{sec_isotoped_rotationmap}
Instead of considering the polynomial $P$, we take the function
$$
g=\sum_{i=0}^{n} z_i^{a_i}-f(r),
$$ 
where $r=\sqrt{\sum_{i=0}^{n}|z_i|^2}$ and the function $f$ is a real valued
function to be specified later. We denote the zero set of $g$ by $\tilde
\Sigma_a$. Note that this set is in general not a Stein manifold. We will, however, show that it is
symplectic for suitable $f$, as one might suspect if $f$ varies slowly. To be more precise, take a vector $X \in T\C^{n+1}|_{g^{-1}(0)}$. The condition that $X$ be tangent to $\tilde \Sigma_a$ is
$$
i_X dg=i_X \left( \sum_{i=0}^{n}a_i z_i^{a_i-1}dz_i-\frac{1}{2}\frac{\partial
  f}{\partial r}\sum_{i=0}^{n}(\frac{\bar z_i}{r}dz_i+\frac{z_i}{r}d\bar
z_i) \right) =0
.
$$
Let now $\omega_0$ denote the standard symplectic form on $\C^{n+1}$
and suppose that $\omega_0|_{\tilde \Sigma_a}$ is degenerate for the vector $X$ at some
point of $\tilde \Sigma_a$. Then we have
$$
i_X \omega_0=(\lambda dg +\bar \lambda d \bar g)
$$
for some $\lambda \in \mathbb C$, because we know $\omega_0$ is
non-degenerate on $\C^{n+1}$. Using this relation, we deduce that
$$
i_X d z_j= \frac{2}{i}\left( - \left( \frac{\partial f}{\partial r} \frac{z_j}{2r} \right) (\lambda+\bar \lambda) +a_j \bar
z_j^{a_j-1} \bar \lambda \right).
$$
Now we return to check the tangency condition of $X$. The previous
relations now give us
$$
0=i_X dg= 
\frac{2}{i}
\bar \lambda \left(
\sum_j a_j^2 |z_j|^{2(a_j-1)}-\frac{\partial f}{\partial
  r}\sum_j\frac{a_j}{2r}(z_j^{a_j}+\bar z_j^{a_j})
\right)
$$
The coefficient of $\bar \lambda$ has a term
involving $a_i^2 |z_i|^{2(a_i-1)}$ in it. Now assume that the exponents are
larger than $1$ and that the derivative $\frac{\partial f}{\partial
  r}<1-\epsilon$ for some positive $\epsilon$. This means that the term with
$a_i^2 |z_i|^{2(a_i-1)}$ will dominate for large $r$, i.e. the
coefficient of $\bar \lambda$ will be non-zero and therefore $\bar \lambda=0$.
Since $|\bar \lambda|=|\lambda|$, it follows that $\lambda$ must be zero,
which in turn implies that $X$ is zero. This last step shows that $\tilde \Sigma_a$ can
be made symplectic for suitable $f$. To be more precise we choose $f$ with
the following properties.
\begin{itemize}
\item[1.] The function $f$ is constant $1$ for $r\leq R_0$, where $R_0$ is chosen in such
a way that the above mentioned term will indeed dominate.
\item[2.] For $r\geq R_1>R_0+1$, the function $f$ is constant $0$. Note that this condition
  is not necessary for symplecticity. It will, however, be useful to
  make the rotation maps isotopic to the identity for large radii.
\item[3.] Between $R_0$ and $R_1$, the function $f$ goes smoothly from $1$
  to $0$, connecting smoothly to the already described parts of $f$. We will choose $f$
such that its derivative is smaller than $1-\epsilon$.
\end{itemize}
Now that we know that $\tilde \Sigma_a$ is symplectic, we want to see that the corresponding
rotation map can be isotoped to the identity. First define the map $\phi:
\C^{n+1} \to \C^{n+1}$, sending $(z_0,\ldots,z_n) \mapsto
(\zeta_{a_0}z_0,z_1,\ldots,z_n)$. Now choose the following Hamiltonian function on $\C^{n+1}$;
$$
H=\sum_{i=0}^{n}\frac{\pi}{a_i}|z_i|^2.
$$
The time $t$ flow of the Hamiltonian vector field associated to $H$ induces
the map 
$$
\psi_t: 
(z_0,\ldots,z_n)
\mapsto
(e^{2 \pi i \frac{t}{a_0} } z_0,\ldots,e^{2 \pi i \frac{t}{a_n} } z_n).
$$
Note that this map sends $\tilde \Sigma_a$ to $\tilde \Sigma_a$ for $r>R_1$. Choose a function $h$
that is constant 0 for $0\leq r\leq R_1$ and that increases to $1$ at $r=R_2>R_1$,
after which it is constant $1$. Let $\tilde \psi_t$ denote the time $t$
flow of the Hamiltonian vector field associated to $\tilde H=hH$. The map $\tilde
\psi_t$ sends $\tilde \Sigma_a $ to $\tilde \Sigma_a$ for all radii. By choosing $t_0 \in \Z$
such that $t_0=-1 ~\mathrm{mod}~a_0$ and $t_0=0 ~\mathrm{mod}~a_i$ for $i=1,\ldots,n$, we undo the
rotation in the first coordinate for large radii and hence we see that $\tilde \psi_t$
is the identity near the boundary. Note this choice is not always possible,
but if $a_0$ is relatively prime to $a_i$ for $i=1,\ldots,n$, it
is. Altogether, we have the map 
$$
\tilde \phi=\tilde \psi_{t_0} \circ \phi: ~\tilde \Sigma_a \to \tilde \Sigma_a,
$$
which is the identity near the boundary of $\tilde \Sigma_a$. Also note
that the choice of $t_0$ is not unique.

\subsubsection{Homomorphism on homology induced by the rotation map}
\label{mappingtorusBrieskorn}
We shall take this isotoped rotation map as the monodromy for the page $\tilde \Sigma_a$. In order to invoke Barden's classification result, we need to know what map the monodromy induces on the homology of $\tilde \Sigma_a$. The Wang sequence we discussed in Section \ref{sec_general} gives the homology of the mapping torus. 

First, we observe that $\phi$ and $\tilde \phi$ are isotopic, so they induce the same maps on homology. And we may, in fact, work with the non-deformed Stein manifold $\Sigma_a$ and the rotation map
defined there (which we will also refer to as $\phi$), because $\tilde \Sigma_a$ and
$\Sigma_a$ coincide in ball of radius $R_0$ around the origin as subsets of
$\C^{n+1}$.

These Stein manifolds $\Sigma_a$ have been studied carefully in the past (see
for instance \cite{Hirzebruch}) and many results about their properties,
including their homology, are known. We will give a short summary of some
of the results that we will use. The results that we are listing are from
Hirzebruch-Mayer, \cite{Hirzebruch}, but date back to Pham, see \cite{Pham}.

In the following we will use the group action on $\Sigma_a$ induced by multiplication by roots of unity. To that end, we introduce some notation. The group of $a_j^\text{th}$-roots of unity will be written as
$G_{a_j}\cong \Z_{a_j}$ when we consider it as an abstract group, and we will denote a
generator of $G_{a_j}$ by $w_j$. As a subgroup of $\C^{*}$, we shall write
$\tilde G_{a_j}$. The roots of unity will be indicated by $\zeta_j$. We
will write $G_a=G_{a_1}\oplus G_{a_2} \oplus \cdots \oplus G_{a_n}$. Let us now consider the deformation retract of the Stein manifolds $\Sigma_a$ indicated in the following theorem.
\begin{theorem}
[Pham, see \cite{Hirzebruch} and \cite{Pham}] The set $U_a=\{ z \in \Sigma_a ~|~z_j^{a_j}\geq 0
\text{ for all } j \}$ is a deformation retract of $\Sigma_a$. This deformation is
compatible with the group action mentioned above.
\end{theorem}
We can parametrize the set $U_a$ in the following way,
$$
U_a=\{ (\zeta_0 t_0, \ldots,\zeta_n t_n) \in \C^{n+1}|~\zeta_j \in \tilde G_{a_j},~t_j\geq 0
\text{ and }\sum_{i=0}^nt_j^{a_j}=1\}.
$$
On the other hand, note that the join $G_{a_0}*\cdots*G_{a_n}$ may be written as 
$$
\tilde G_{a_0}*\cdots*\tilde G_{a_n}=\{(\zeta_0 t_0, \ldots,\zeta_n t_n) \in \C^{n+1}~|~\zeta_j \in
\tilde G_{a_j},~t_j\geq 0
\text{ and }\sum_{i=0}^nt_j=1  \}.
$$
These sets can be identified if we rescale the $t_j$'s. Notice that this
identification is compatible with the group action, because $G_a$ acts only
on the roots of unity.

General theory gives us that the join $G_{a_0}*\cdots*G_{a_n}$ is an $n$-dimensional simplicial complex with
an $n$-simplex for each element in $G_a$. This is again
compatible with the group action in the following sense. Let $e$ denote the
simplex corresponding to $1\in G_a$. The other simplices are obtained by
letting $G_a$ act. In other words, the simplicial chain complex in degree
$n$ can be written as
$$
C_n(U_a)=\Z(G_a) e,
$$
where $\Z(G_a)$ denotes the group ring of $G_a$.

Now define
$$
h:=(1-w_0)(1-w_1)\ldots(1-w_n)e.
$$
In the lecture notes of Hirzebruch and Mayer \cite{Hirzebruch} it is shown that $h$ is a cycle. In fact, one can establish an isomorphism (see \cite{Hirzebruch} for more details)
$$
\tilde H_n(U_a) \cong \Z(G_a)h
$$
coming from the homomorphism
\begin{eqnarray*}
\Phi:~C_n(U_a) \cong \Z(G_a) & \to & \Z(G_a) h \\
w & \mapsto & w h.
\end{eqnarray*}
The kernel of $\Phi$ is the ideal $I_a$ generated by
$$
1+w_j+w_j^2+\cdots+w_j^{a_j-1} \text{ for }j=0,\ldots,n.
$$
Let us consider the basis of $\tilde H_n(U_a)$ represented by elements in $C_n(U_a)$ of the form
\begin{equation}
\label{eqn_basis}
w_0^{k_0}w_1^{k_1} \cdot \ldots \cdot w_n^{k_n} \text{ with } 0\leq k_j \leq
a_j-2 \text{ for }j=0,\ldots,n.
\end{equation}
With respect to this basis, we can give a matrix representation of $\phi_*$, the isomorphism on homology induced by the rotation map $\phi$. The "rotation" map $\phi$ corresponds to multiplication by $w_0$ on $C_n(U_a)$. That is to say that $\phi$ shifts the basis in Formula \ref{eqn_basis} by $w_0$. For the induced map on homology, we use the ideal $I_a$ to simplify the results if necessary, for instance
\begin{eqnarray*}
w_0 & \mapsto & w_0^2 \\
w_0^{a_0-2} & \mapsto & w_0^{a_0-1}\equiv -1-w_0-\ldots -w_0^{a_0-2} \text{ mod }I_a.
\end{eqnarray*}
Hence the matrix representation of $\phi_{*}$ consists of $(a_1-1)\cdot
\ldots \cdot (a_n-1)$ blocks on the diagonal that look like the $(a_0-1)\times (a_0-1)$-matrix
$$
\left(
\begin{array}{ccccc}
 0 & 0 & \cdots & 0 & -1\\
 1 & 0 & \cdots & 0 & ~ \\
 0 & 1 & \ddots & \vdots & \vdots \\
 \vdots & \ddots & \ddots  & 0 & ~ \\
 0 & \cdots & 0 & 1 & -1 \\ 
\end{array}
\right)
$$
if we order the basis by its degree in $w_1$, then by its degree in $w_2$
and so on.

The above representation of $\phi_{*}$ can be used to compute the homology
of the mapping torus
$$
A':=\Sigma_a\times I /\sim,\text{ where }(x,0)\sim(\phi(x),1).
$$
This is done using the Wang sequence. We use the
facts that $H_3(\Sigma_a)=0$ and that $\pi_1(\Sigma_a)=0$ (and hence also
$H_1(\Sigma_a)=0$). The piece that is relevant to us looks like
$$
\label{Wangsequence_mappingtorus}
0 \to H_3(A') \to H_2(\Sigma_a) \stackrel{\phi_{*}-id}{\to} H_2(\Sigma_a)
\to H_2(A') \to 0.
$$
Using the above matrix representation of $\phi_{*}$ we see that
$\phi_{*}-id$ is injective, because the determinant of the associated
matrix is non-zero. Hence we conclude that $H_3(A')=0$ and that
$H_2(A')\cong \mathrm{coker}(\phi_{*}-id)$. We have 
$$
H_2(A')\cong \mathrm{coker}(\phi_{*}-id) \cong 
\underset{(a_1-1) \cdot \ldots \cdot (a_n-1)\text{times}}
{\underbrace{\Z_{a_0} \oplus \cdots \oplus \Z_{a_0}. }}
$$
Indeed, each block of the matrix representation of $\phi_{*}-id$
corresponding to the above block has a cokernel isomorphic to $\Z_{a_0}$,
which can be seen by performing Gauss elimination. Together with the
discussion at the beginning of this section this gives us the homology of
the mapping torus of $\tilde \Sigma_a$ with monodromy $\tilde \phi$. Let $A$ denote this
mapping torus,
$$
A=\tilde \Sigma_a \times I /\sim,\text{ where }(x,0)\sim(\tilde\phi(x),1).
$$
Then we have 
\begin{equation}
\label{homologyBrieskornmappingtorus}
H_2(A)\cong \underset{(a_1-1) \cdot \ldots \cdot (a_n-1)\text{times}}
{\underbrace{\Z_{a_0} \oplus \cdots \oplus \Z_{a_0} }}.
\end{equation}
The homotopy exact sequence of the fibration $A \to S^1$ shows that
$\pi_1(A)\cong \Z$, so we see that $H_1(A)\cong \Z$ as well. All higher
homology groups (grade larger than two) are zero.

\subsubsection{Homology of the open book}
Now we choose suitable exponents for the Brieskorn varieties and  use them to give the remaining prime manifolds contact open books.

First of all, we consider the Brieskorn variety $\tilde \Sigma_a$ with exponents
$a_0=p^k$, $a_1=3$ and $a_2=2$, where $p$ is a prime different from $2$ and
$3$, and $k$ some positive integer. Notice that the associated Brieskorn
manifold $K$ is then a homology sphere, i.e.~$H_1(K)=0$. The set $A$ denotes the
mapping torus of $\tilde \Sigma_a$ with monodromy $\tilde\phi$ as in
the previous section. As is our convention, we define $B:=K\times D^2$ and set
$X:=A\cup_\partial B$.

The arguments from Section \ref{sec_general} show that $X$ is simply connected. By Poincar\'e duality we see that
$H_4(X)=0$, and since $K$ is a homology sphere we also have $H_2(A\cap
B)=0$. Consider the following piece of the Mayer-Vietoris sequence,
$$
0 \to \underset{\cong \Z_{p^k}\oplus
\Z_{p^k}}{H_2(A)}\oplus \underset{\cong 0}{H_2(B)} \to H_2(X) \to 0.
$$
Here we have used the arguments from Section
\ref{sec_general} to split off a part of the sequence. We
see directly that $H_2(X)\cong \Z_{p^k}\oplus \Z_{p^k}$. In particular, the rank of $H_2(X)$ is zero, so $H_3(X)=0$ as
well by Poincar\'e duality and the universal coefficient theorem. This shows that the prime
manifolds $M$ with $H_2(M)\cong \Z_{p^k}\oplus \Z_{p^k}$ admit contact open
books for $p\neq 2,3$. The binding is a Brieskorn homology sphere of the
form $\Sigma(p^k,3,2)$, and the page is the Brieskorn variety $\tilde
\Sigma_a$. Together with our earlier results, this covers all prime manifolds
except those with $2$- or $3$-torsion in their second homology group. To get them,
we consider Brieskorn varieties with different exponents.

First we shall tackle the case of $2$-torsion in homology. Consider the
Brieskorn varieties $\tilde \Sigma_a$ with exponents $a_0=2^k$, $a_1=3$ and
$a_2=3$. Since the exponents are not relatively prime, we cannot conclude
that $K$ is a homology sphere. We can, however, compute the homology of $K$
by using the algorithm of Randell \cite{Randell}. We get $H_1(K)\cong \Z_{2^k}\oplus \Z_{2^k}$. 

Let $X:=A\cup_\partial B$ be the open book as before, but now with the new exponents. If we consider the Mayer-Vietoris sequence for $(A,B)$ in $X$ with rational coefficients, we easily see that the rank of $H_3(X)$ is
zero. Together with the arguments from Section
\ref{sec_general} this reduces the remaining part of the
Mayer-Vietoris sequence for $(A,B)$ with integer coefficients to
$$
0 \to \underset{\cong \Z_{2^k}^2}{H_2(A\cap B)} \overset{i\oplus j}\to
\underset{\cong \Z_{2^k}^4 }{H_2(A)}\oplus
\underset{=0}{H_2(B)} \to H_2(X) \to 0.
$$
We have used to K\"unneth formula to determine $H_2(A\cap B)$. The rank of
$H_1(K)$ is zero, so by Poincar\'e duality $H_2(K)=0$ and hence we also have
$H_2(B)=0$. Formula (\ref{homologyBrieskornmappingtorus}) gives the homology
of $A$. Injectivity of the map ${i\oplus j}$ means that we can represent this map by a $(4\times 2)$ matrix which has a $(2\times 2)$ subdeterminant that is a unit in $\Z_{2^k}$. Hence we see that we can extend this matrix to form a basis of $\Z_{2^k}^4$. So after applying a basis transformation on $\Z_{2^k}^4$, we see that
$$
\im {i\oplus j}=\Z_{2^k}\times \Z_{2^k}\times \{ 0\} \times \{ 0 \}.
$$
Hence by exactness, we obtain $H_2(X)\cong \Z_{2^k}^2$.

The arguments for the $3$-torsion case are almost completely the same. The
exponents for $\Sigma_a$ are different, of course. We take $a_0=3^k$, $a_1=4$
and $a_2=2$. As before we use the algorithm of Randell \cite{Randell} to
compute the homology of the Brieskorn manifold $K$. This time we get
$H_1(K)\cong \Z_{3^k}$. Formula (\ref{homologyBrieskornmappingtorus}) shows
that $H_2(A)=\Z_{3^k}^3$. Again, using the arguments from Section
\ref{sec_general} we can split off a part of the
Mayer-Vietoris sequence. By tensoring with $\Q$ we see that the rank of $H_2(X)$
is zero, and hence $H_3(X)=0$. This reduces the sequence to
$$
0 \to \underset{\cong \Z_{3^k}}{H_2(A\cap B)} \stackrel{i\oplus j}{\to} \underset{\cong
  \Z_{3^k}^3}{H_2(A)} \oplus \underset{\cong 0}{H_2(B)} \to H_2(X) \to 0.
$$
The map $i\oplus j$ is injective, so $i\oplus j(\bar 1)=(a,b,c)$ is an element of order $3^k$. This means that one of the elements $a,b,c$ is a unit in $\Z_{3^k}$. Therefore we can include the vector $(a,b,c)$ into a basis of $\Z_{3^k}^3$. With respect to this basis we have $i\oplus j(\bar 1)=(1,0,0)$. By exactness, we see directly that $H_2(X)\cong \Z_{3^k}\oplus \Z_{3^k}$. 

\begin{remark}
An easier way to see that these prime manifolds admit contact structures is by considering Brieskorn manifolds. Namely, we have
$$
H_2(\Sigma(p^k,3,3,3))\cong \Z_{p^k}\oplus \Z_{p^k} \text{ for } p\text{ not divisible by }3
$$
and
$$
H_2(\Sigma(p^k,2,4,4))\cong \Z_{p^k}\oplus \Z_{p^k} \text{ for } p\text{ not divisible by }2.
$$
This can be shown by applying Randell's algorithm \cite{Randell}. Of course, we do not obtain the abstract open books in this way.
\end{remark}

\section{Conclusion and discussion}
In Section \ref{sec_openbook_free} and Section \ref{sec_openbook_torsion} we constructed abstract contact open books for the prime manifolds in Barden's classification. Note that we can easily obtain an abstract contact open book for $S^5$. Simply take $D^4$ with standard symplectic structure as page and use the identity as monodromy. In view of Remark \ref{rmk_booksum}, this gives abstract contact open books for all simply-connected five-manifolds that admit an almost contact structure. Moreover, we can realize a contact structure with any admissible Chern class, since a non-zero Chern class can only come from an $S^2\times S^3$- or an $S^2\tilde \times S^3$-factor. For the latter two manifolds we have shown that we can realize all possible Chern classes. In \cite{Geiges_five} Lemma~7, it is shown that, for an oriented five-manifold, the almost contact structure is completely determined by the first Chern class. 

We can change the orientation by replacing a contact form $\alpha$ by $-\alpha$. Hence we get a contact open book for every homotopy class of almost contact structures on a simply-connected five-manifold. This completes our alternative proof of Theorem \ref{thm_main}.

\begin{remark}
In our construction there is still a lot of freedom left, even though we took very explicit cases. For $S^2\times S^3$ and $S^2\tilde \times S^3$ we can, for instance, vary the page of the abstract open book but keep the Chern class fixed. This can for instance be done by adding two stabilizations to the Legendrian unknot used for the handle attachment; by adding one stabilization on the left side and one on the right side of the Legendrian unknot, we fix the rotation number, but decrease the framing ($tb-1$) by $2$. The resulting abstract contact open books have the same Chern class, but are they contactomorphic?

For the other prime manifolds, we can vary the monodromy in the following way. The parameter $t_0$ we used in isotoping the "rotation" map to the identity near the boundary in Section \ref{sec_isotoped_rotationmap} is not unique. The obvious question is, whether different choices can lead to different contact structures on the same manifold. Here one should note that although we used a Hamilton vector field for the isotopy, we did not use one with compact support. Hence we could get different maps that are not symplectically isotopic relative to boundary. 
\end{remark}

\bibliographystyle{cdraifplain}

\begin{thebibliography}{10}
\bibitem{ACampo}
N.~A'Campo, \emph{Feuilletages de codimension $1$ sur des variétés de dimension $5$}, C. R. Acad. Sci. Paris Sér. A-B 273 (1971), A603--A604.

\bibitem{Barden}
D.~Barden, \emph{Simply connected five-manifolds}, Ann. of Math. (2) 82
(1965), 365--385.

\bibitem{Geiges_five}
H.~Geiges, \emph{Contact structures on $1$-connected $5$-manifolds}, Mathematika 38 (1991), no. 2, 303--311. 

\bibitem{Giroux_talk}
E.~Giroux, J.~Mohsen, \emph{Contact structures and symplectic fibrations
  over the circle}, lecture notes.

\bibitem{Gompf_Stein}
R.~Gompf, \emph{Handlebody construction of Stein surfaces}, Ann. of Math. (2) 148 (1998), no. 2, 619--693.

\bibitem{Gompf_Stipsicz}
R.~Gompf, A.~Stipsicz, \emph{$4$-manifolds and Kirby
calculus}, Graduate  Studies in Mathematics, 20. American Mathematical Society, Providence, RI, 1999.

\bibitem{Hirzebruch}
F.~Hirzebruch, K.~Mayer, \emph{{${\rm O}(n)$}-{M}annigfaltigkeiten,
  exotische {S}ph\"aren und {S}ingularit\"aten}, Lecture Notes in Mathematics,
  No. 57, Springer-Verlag, Berlin, 1968.

\bibitem{Pham}
F.~Pham, \emph{Formules de Picard-Lefschetz généralisées et ramification
  des intégrales}, Bull. Soc. Math. France 93 (1965) 333--367.

\bibitem{Randell}
R.~Randell, \emph{The homology of generalized Brieskorn manifolds}, Topology 14 (1975), no. 4, 347--355.

\bibitem{Thurston_Winkelnkemper}
W.~Thurston, H.~Winkelnkemper, \emph{On the existence of contact forms}, Proc. Amer. Math. Soc. 52 (1975), 345--347.

\bibitem{Wang}
H-C.~Wang, \emph{Homology of fibre bundles Duke}, Math Journal 16, 33-38.

\end{thebibliography}

\end{document}